\documentclass{amsart}

\usepackage{mathrsfs}
\usepackage{stmaryrd}
\usepackage{enumerate}
\usepackage{tikz-cd}
\usepackage{tikz}
\usetikzlibrary{calc, shapes, backgrounds,arrows,positioning,plotmarks}
\tikzset{>=stealth',
  head/.style = {fill = white, text=black},
  plaque/.style = {draw, rectangle, minimum size = 10mm}, 
  pil/.style={->,thick},
  junct/.style = {draw,circle,inner sep=0.5pt,outer sep=0pt, fill=black}
  }
\usepackage[all]{xy}
\usepackage{aliascnt}

\usepackage{fullpage}
\usepackage{verbatim}
\usepackage{amssymb}

\newcommand{\arxiv}[1]{\href{https://arxiv.org/abs/#1}{\texttt{arXiv:#1}}}

\definecolor{darkblue}{rgb}{0.0,0,0.7}

\newcommand{\newword}[1]{\textcolor{darkblue}{\textbf{\emph{#1}}}}

\usepackage[colorlinks=true, linkcolor=blue, citecolor=darkblue]{hyperref}

\newtheorem{theorem}{Theorem}

\newaliascnt{lemma}{theorem}

\aliascntresetthe{lemma}

\newaliascnt{corollary}{theorem}

\aliascntresetthe{corollary}

\newaliascnt{proposition}{theorem}

\aliascntresetthe{proposition}

\newaliascnt{conjecture}{theorem}

\aliascntresetthe{conjecture}

\theoremstyle{definition}

\newaliascnt{definition}{theorem}

\aliascntresetthe{definition}

\newaliascnt{remark}{theorem}

\aliascntresetthe{remark}

\newaliascnt{example}{theorem}

\aliascntresetthe{example}

\newaliascnt{notation}{theorem}

\aliascntresetthe{notation}

\newcommand{\Z}{\mathbb{Z}}

\newcommand{\C}{\mathbb{C}}

\newcommand{\bP}{\mathbb{P}}

\newcommand{\cE}{\mathcal{E}}
\newcommand{\cO}{\mathcal{O}}

\newcommand{\cP}{\mathcal{P}}

\newcommand{\cQ}{\mathcal{Q}}
\newcommand{\cL}{\mathcal{L}}

\newcommand{\mfG}{\mathfrak{G}}
\newcommand{\mfS}{\mathfrak{S}}

\newcommand{\xx}{\mathbf{x}}

\DeclareMathOperator{\supp}{supp}

\DeclareMathOperator{\Fl}{\mathcal{F}\ell}

\newif\ifhascomments \hascommentstrue
\ifhascomments
  \newcommand{\matt}[1]{{\color{red}[[\ensuremath{\spadesuit\spadesuit\spadesuit} #1]]}}
  \newcommand{\oliver}[1]{{\color{blue}[[\ensuremath{\clubsuit\clubsuit\clubsuit} #1]]}}
\else
  \newcommand{\oliver}[1]{}
  \newcommand{\matt}[1]{}
\fi

\begin{document}

\title{Proof of a conjectured M\"obius inversion formula for Grothendieck polynomials}

\date{\today}

\author{Oliver Pechenik}
\thanks{OP was partially supported by a Discovery Grant (RGPIN-2021-02391) and Launch Supplement (DGECR-2021-00010) from the Natural Sciences and Engineering Research Council of Canada.}
\address[OP]{Department of Combinatorics \& Optimization, University of Waterloo, Waterloo ON N2L3G1, Canada}
\email{oliver.pechenik@uwaterloo.ca}

\author{Matthew~Satriano}
\thanks{MS was partially supported by a Discovery Grant from the
  National Science and Engineering Research Council of Canada and a Mathematics Faculty Research Chair.}
\address[MS]{Department of Pure Mathematics, University
  of Waterloo, Waterloo ON N2L3G1, Canada}
\email{msatrian@uwaterloo.ca}


\makeatletter
\@namedef{subjclassname@2020}{%
  \textup{2020} Mathematics Subject Classification}
\makeatother


\begin{abstract}
	Schubert polynomials $\mfS_w$ are polynomial representatives for cohomology classes of Schubert varieties in a complete flag variety, while Grothendieck polynomials $\mfG_w$ are analogous representatives for the $K$-theory classes of the structure sheaves of Schubert varieties. In the special case that $\mfS_w$ is a multiplicity-free sum of monomials, K.~M\'esz\'aros, L.~Setiabrata, and A.~St.~Dizier conjectured that $\mfG_w$ can be easily computed from $\mfS_w$ via M\"obius inversion on a certain poset. We prove this conjecture. Our approach is to realize monomials as Chow classes on a product of projective spaces and invoke a result of M.~Brion on flat degenerations of such classes.
\end{abstract}

\maketitle

\section{Introduction}\label{sec:intro}

The \newword{flag variety} $\Fl_n$ is the moduli space of complete flags
\[
V_0 \subset V_1 \subset \dots \subset V_n = \C^n
\]
of nested vector subspaces of $\C^n$, where $\dim V_i = i$. The flag variety is stratified by its \emph{Schubert varieties} $X_w$ for $w \in S_n$. Taking Poincar\'e duals of these subvarieties yields a distinguished ``Schubert basis'' of the integral cohomology ring $H^\star(\Fl_n)$. Similarly, taking classes of their structure sheaves gives a distinguished Schubert basis of the $K$-theory ring $K^0(\Fl_n)$.

It is an important problem to be able to compute in these cohomology and $K$-theory rings. In particular, one would like a positive combinatorial formula for the structure coefficients of these rings with respect to their distinguished bases. For cohomological Schubert classes pulled back from a complex Grassmannian, the famous Littlewood--Richardson rule provides such a formula. Analogous formulas are known in certain other special cases, but the general problem remains very open, despite much attention.

Towards solving this problem, A.~Lascoux and M.-P.~Sch\"utzenberger \cite{Lascoux.Schutzenberger:Schubert,Lascoux.Schutzenberger} introduced polynomial representatives, called \emph{Schubert} and \emph{Grothendieck polynomials}, for the Schubert classes in cohomology and $K$-theory, respectively. These are representatives in the sense that the Schubert polynomials $\mfS_w \in \Z[x_1, \dots, x_n]$ satisfy
\[
\mfS_u \cdot \mfS_v = \sum_w c_{u,v}^w \mfS_w
\]
if and only if the cohomological Schubert classes $\sigma_w \in H^\star(\Fl_n)$ (for sufficiently large $n$) satisfy
\[
\sigma_u \cdot \sigma_v = \sum_w c_{u,v}^w \sigma_w,
\] with an analogous statement for $K$-theoretic classes and Grothendieck polynomials $\mfG_w \in \Z[x_1 \dots, x_n]$.

Beautiful combinatorial formulas for Schubert and Grothendieck polynomials are known (see, e.g., \cite{Knutson.Miller,Weigandt}. However, much of their structure remains mysterious, especially in the case of Grothendieck polynomials, despite recent work such as \cite{Pechenik.Speyer.Weigandt,Hafner}. For background on Schubert and Grothendieck polynomials and associated geometry, see, e.g., \cite{Fulton97,Manivel,Knutson.Miller}.

In certain cases (classified in \cite{Fink.Meszaros.StDizier} via permutation pattern avoidance), the Schubert polynomial $\mfS_w$ is a multiplicity-free sum of monomials (i.e., all nonzero coefficients are $1$). In this setting, K.~M\'{e}sz\'{a}ros, L.~Setiabrata, and A.~St.~Dizier \cite{MSSD} conjectured that the monomial expansion of the Grothendieck polynomial $\mfG_w$ is easily determined by $\mfS_w$ via M\"obius inversion. Our main result is to prove their conjecture.

More specifically, for $w \in S_n$, they introduce an integer vector $\textrm{wt}(\overline{D(w)})$ (\cite[\S3]{MSSD}) and a poset (\cite[Definition 4.1]{MSSD})
\[
P_w = \{\beta\in\Z^n\mid\alpha\leq\beta\leq\textrm{wt}(\overline{D(w)})\textrm{\ for\ some\ }\alpha\in\supp(\mfG_w)\} \; \sqcup \; \{ \hat{0} \},
\] where $(\alpha_1, \dots \alpha_n)\in\supp(\mfG_w)$ denotes that the monomial $x_1^{\alpha_1} x_2^{\alpha_2}\cdots x_n^{\alpha_n}$ has nonzero coefficient in the Grothendieck polynomial $\mfG_w$. Here, the poset structure on integer vectors in $P_w$ is componentwise comparison, while $\widehat{0}$ is the unique minimum element.
M\'esz\'aros--Setiabrata--St.~Dizier then conjecture that, for $w$ such that $\mfS_w$ is multiplicity-free, the corresponding Grothendieck polynomial $\mfG_w$ is computed via M\"obius function on the poset $P_w$. The following theorem proves their conjecture.

\begin{theorem}[{\cite[Conjecture 1.5]{MSSD}}]
\label{thm:main}
Let $w$ be a permutation such that all nonzero coefficients of $\mfS_w$ equal $1$. If $\mu_w$ is the M\"obius function of $P_w$, then
\[
\mfG_w=-\sum_{\beta\in P_w\smallsetminus\widehat{0}}\mu_w(\widehat{0},\beta)x_1^{\beta_1} x_2^{\beta_2} \cdots x_n^{\beta_n}.
\]
\end{theorem}

Our proof of Theorem~\ref{thm:main} appears in Section~\ref{sec:proof}. Our approach is to realize the monomials of $\mathfrak{S}_w$ as Chow classes on a product of projective spaces. From there, our main tools are results of M.~Brion \cite{Brion} and A.~Knutson \cite{Allen} on subvarieties of generalized flag varieties whose Chow classes are multiplicity-free in the Schubert basis.
For examples of the application of Theorem~\ref{thm:main}, see Section~\ref{sec:example} or \cite[Example~4.2 and Figure~7]{MSSD}.

\section{Proof of \autoref{thm:main}}\label{sec:proof}

Fix $n\geq 1$ and $w\in S_n$. For any $m\geq n$, let $\Fl(n,m)$ be the partial flag variety parameterizing complete flags of an $n$-dimensional vector space embedded in $\C^m$. The partial flag variety $\Fl(n,m)$ comes equipped with a universal flag
\[
0=\cE_0\subset \cE_1\subset \cE_2\subset\dots\subset \cE_n \subset \cO_{\Fl(n,m)}^{\oplus m}.
\]
For $0<i\leq n$, let $\cQ_i:=\cE_i/\cE_{i-1}$ be the $i$th quotient line bundle. The Schubert polynomial $\mfS_w\in\Z[x_1,\dots,x_n]$ has the property that
\[
\mfS_w(c_1(\cQ_1^\vee),\dots,c_1(\cQ_n^\vee))=\sigma_w \in H^{2\ell(w)}(\Fl(n,m))
\]
where $\sigma_w$ is the class of the Schubert variety $X_w$ and $\ell(w)$ denotes the Coxeter length of the permutation $w$ (see, e.g., \cite[p.~181]{Fulton97}). Similarly, the Grothendieck polynomial $\mfG_w\in\Z[x_1,\dots,x_n]$ has the property that
\[
\mfG_w(c_1^K(\cQ_1),\dots,c_1^K(\cQ_n))=[\cO_{X_w}]\in K^0(\Fl(n,m))
\]
where $c_1^K$ denotes the $K$-theoretic Chern class (see, e.g., \cite[p.~72]{Buch}). Here, for a line bundle $\cL$ over a base $X$, the $K$-theoretic Chern class is defined by $c_1^K(\cL) = 1 - [\cL] \in K^0(X)$.

Let $d\geq\deg\mfG_w$ and consider the $n$-fold product $(\bP^d)^n$. (Throughout, we write $\bP^k$ to mean specifically $\C\bP^k$.) Let $\pi_i\colon (\bP^d)^n\to\bP^d$ be the projection to the $i$th factor and let
\[
\cL_i:=\pi_i^*\cO_{\bP^d}(-1)\subset\cO_{(\bP^d)^n}^{d+1}.
\]
Then $\bigoplus_{i=1}^n\cL_i\subset \cO_{(\bP^d)^n}^{n(d+1)}$. Moreover, for all $m\geq n(d+1)$, 
we have a flag 
\[
\cL_1\subset \cL_1\oplus\cL_2\subset\dots\subset \bigoplus_{i=1}^n\cL_i\subset \cO_{(\bP^d)^n}^m.
\]
This defines a map $f\colon(\bP^d)^n\to\Fl(n,m)$.
By construction, 
$f^*(\cQ_i)=\cL_i$.
Consider the maps
\[
f^*\colon H^\star(\Fl(n,m))\to H^\star((\bP^d)^n)\simeq \Z[y_1,\dots,y_n]/(y_1^{d+1},\dots,y_n^{d+1})
\]
and 
\[
f_K^*\colon K^0(\Fl(n,m))\to K^0((\bP^d)^n)\simeq \Z[z_1,\dots,z_n]/((z_1-1)^{d+1},\dots,(z_n-1)^{d+1})
\]
on cohomology and $K$-theory; here $y_i$ is the class of $[\bP^n\times\dots\times\bP^n\times\bP^{n-1}\times\bP^n\times\dots\times\bP^n]$ with $\bP^{n-1}$ in the $i$th factor and $z_i=[\pi_i^*\cO_{\bP^n}(-1)]=[\cL_i]$. Note that
\[
f^*(-c_1(\cQ_i))=c_1(\cL_i^\vee)=c_1(\pi_i^*\cO(1))=y_i
\]
and that 
\[
f_K^*(c_1^K(\cQ_i))=c_1^K(\cL_i)=1-[\cL_i]=1-z_i.
\]
It follows that
\begin{equation}\label{eqn:pullbackK-class}
f^*\sigma_w=\mfS_w(y_1,\dots,y_n)\quad\textrm{and}\quad
f_K^*[\cO_{X_w}]=\mfG_w(1-z_1,\dots,1-z_n).
\end{equation}
Furthermore, since $d\geq\deg\mfG_w\geq\deg\mfS_w$, no term of $\mfS_w\in\Z[x_1,\dots,x_n]$ (respectively, $\mfG_w\in\Z[x_1,\dots,x_n]$) is killed when substituting $y_i$ (respectively, $1-z_i$) in place of $x_i$.

Notice that each monomial $y_1^{a_1}\dots y_n^{a_n}$ is a Schubert class in the homogeneous space $(\bP^d)^n$, namely the class $[\bP^{d-a_1}\times\dots\times\bP^{d-a_n}]$. Thus, under the hypotheses of Theorem~\ref{thm:main}, $f^*\sigma_w$ is a multiplicity-free sum of Schubert classes in $H^\star((\bP^d)^n)$. Since the cohomology and Chow rings are isomorphic for $(\bP^d)^n$, with the isomorphism matching up Schubert classes, the Chow class of $f^{-1}X_w$ is also a multiplicity-free sum of Schubert classes. M.~Brion \cite{Brion} showed that such classes in flag varieties have nice flat degenerations establishing various very special properties.  Building on this work, A.~Knutson \cite[Theorem 3]{Allen}  explained how one may compute the corresponding $K$-class from the multiplicity-free Chow expansion in this setting. In our case, Knutson's formula gives
\begin{equation}\label{eqn:pullbackK-class-Mobius}
f_K^*[\cO_{X_w}]=\sum_{W\in\cP}\mu_{\cP}(W)[\cO_W];
\end{equation}
here $\cP$ is the poset of subvarieties $\bP^{b_1}\times\dots\times\bP^{b_n}$, ordered by inclusion, such that there exists $\bP^{c_1}\times\dots\times\bP^{c_n}$ in the support of $f^*[X_w]$ with $b_i\leq c_i$ for all $i$. Knutson uses the convention that for all $W\in\cP$, we have $\sum_{Z\geq W}\mu_{\cP}(Z)=1$. (Technically, Knutson's formula is stated for $K_0((\bP^d)^n)$ rather than $K^0((\bP^d)^n)$, but the distinction is immaterial for flag varieties, which are smooth; for further discussion of the $K$-homology versus $K$-cohomology issue, see \cite{Knutson:Osaka}.)

Identifying $\bP^{b_1}\times\dots\times\bP^{b_n}$ with the integer vector $(d-b_1,\dots,d-b_n)$, we see $\cP$ is isomorphic to the poset $\cP'$ of tuples $(b_1,\dots,b_n) \in \Z^n$ such that:
\begin{itemize}
	\item[(i)] all $b_i\leq d$, and 
	\item[(ii)]there exists a monomial $x_1^{a_1}\dots x_n^{a_n}$ in the support of $\mfS_w$ with $a_i\leq b_i$ for all $i$;
\end{itemize}
since the lowest-degree terms of $\mfG_w$ coincide with $\mfS_w$, condition (ii) is equivalent to requiring the existence of a monomial $x_1^{a_1}\dots x_n^{a_n}$ in the support of $\mfG_w$ with $a_i\leq b_i$ for all $i$. Under this isomorphism between $\cP$ and $\cP'$, we must now use the convention that $\sum_{\alpha\geq\beta}\mu_{\cP'}(\beta)=0$ for all $\alpha$.

Iteratively using the short exact sequences
\[
0\to\cO_{\bP^{r}}(-1)\to\cO_{\bP^{r}}\to\cO_{\bP^{r-1}}\to0,
\]
one readily checks that in $K^0(\bP^d)$, we have
\begin{equation}\label{eqn:OPr-Kthy}
[\cO_{\bP^{r}}]=\sum_{i=0}^{d-r}(-1)^i{d-r\choose i}[\cO_{\bP^d}(-i)]=(1-[\cO_{\bP^d}(-1)])^{d-r}.
\end{equation}
Combining equations \eqref{eqn:pullbackK-class}, \eqref{eqn:pullbackK-class-Mobius}, and \eqref{eqn:OPr-Kthy}, we find
\begin{align*}
\mfG_w(1-z_1,\dots,1-z_n)=f_K^*[\cO_{X_w}] &=\sum_{(b_1,\dots,b_n)\in\cP'}\mu_{\cP'}((b_1,\dots,b_n))[\cO_{\bP^{d-b_1}\times\dots \times \bP^{d-b_n}}]\\
&=\sum_{(b_1,\dots,b_n)\in\cP'}\mu_{\cP'}((b_1,\dots,b_n))\prod_{i=1}^n(1-z_i)^{b_i}.
\end{align*}
We see then that the coefficient of the monomial $x^{b_1}\cdots x^{b_n}$ in $\mfG_w(x_1,\dots,x_n)$ is 0 if $(b_1,\dots,b_n)\notin\cP'$ and is $\mu_{\cP'}((b_1,\dots,b_n))$ otherwise. Applying \cite[Theorem 1.2]{MSSD21}, it follows that $\mu_{\cP'}((b_1,\dots,b_n))=0$ unless $(b_1,\dots,b_n)\leq\textrm{wt}(\overline{D(w)})$. Hence, if we let $\cP''$ be the subposet of $\cP'$ where we additionally require $(b_1,\dots,b_n)\in\cP''$ to satisfy $(b_1,\dots,b_n)\leq\textrm{wt}(\overline{D(w)})$, then we see
\[
\mfG_w(x_1,\dots,x_n)=\sum_{(b_1,\dots,b_n)\in\cP''}\mu_{\cP''}((b_1,\dots,b_n))x_1^{b_1}\cdots x_n^{b_n}.
\]
We now simply note that $\cP''=P_w\smallsetminus\widehat{0}$, where $P_w$ is as defined in Section~\ref{sec:intro}. The authors of \cite{MSSD} use the convention $\mu_{P_w}(\widehat{0},\widehat{0})=1$ and $\sum_{\alpha\geq\beta}\mu_{P_w}(\widehat{0},\beta)=0$ for all $\alpha$. It follows that 
\[
\mu_{P_w}(\widehat{0},(b_1,\dots,b_n))=-\mu_{\cP''}((b_1,\dots,b_n)),
\]
and hence
\[
\mfG_w(x_1,\dots,x_n)=-\sum_{\beta\in P_w\smallsetminus\widehat{0}}\mu_{P_w}(\widehat{0},\beta)x^\beta,
\]
which establishes Theorem~\ref{thm:main}. \qed

\section{An example}\label{sec:example}
Given a permutation $w \in S_n$, we write it in one-line notation as the string $w(1)w(2) \ldots w(n)$.
Let $w = 143562 \in S_6$. In this case, the Schubert polynomial $\mathfrak{S}_{143562}$ has a multiplicity-free monomial expansion by \cite[Theorem~1.1]{Fink.Meszaros.StDizier}. Given an integer vector $a = (a_1, a_2, \ldots, a_n)$, we write $\xx^a = x_1^{a_1} x_2^{a_2} \cdots x_n^{a_n}$ for concision. One may calculate explicitly that
\begin{align*}
\mathfrak{S}_{143562} &= \xx^{(2,0,1,1,1)} + \xx^{(1,1,1,1,1)} + \xx^{(2,1,1,1,0)} + \xx^{(1,2,1,1,0)} + \xx^{(2,1,1,0,1)} \\ &+ \xx^{(1,2,1,0,1)} + \xx^{(2,1,0,1,1)} + \xx^{(1,2,0,1,1)} + \xx^{(0,2,1,1,1)}. 
\end{align*}
 It is easily computed that the vector $\gamma = \mathrm{wt}(\overline{D(143562)}) = (2,2,1,1,1)$.
 
 Our poset $P_{143562}$ then has Hasse diagram
 \[\scalebox{0.7}
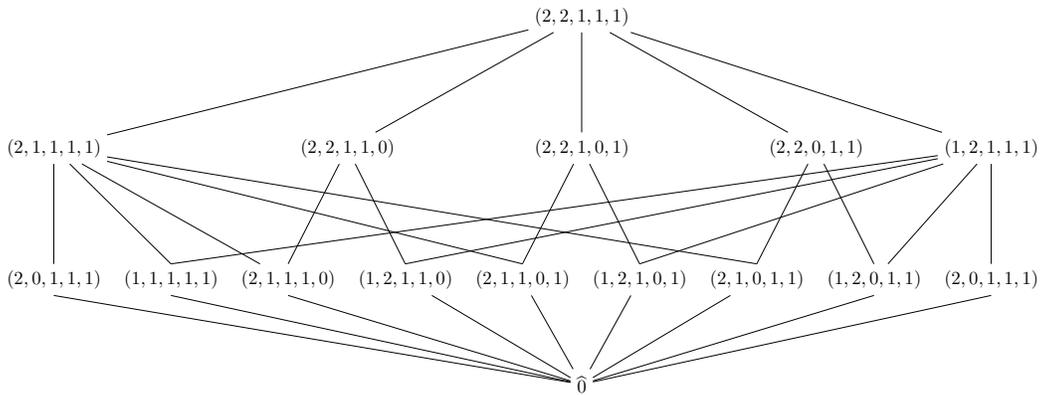
{
 \begin{tikzpicture}
	\node (20111) at (0,0) {$(2,0,1,1,1)$};
	\node[right= 0.2 of 20111] (11111) {$(1,1,1,1,1)$};
	\node[right= 0.2 of 11111] (21110) {$(2,1,1,1,0)$};
	\node[right= 0.2 of 21110] (12110) {$(1,2,1,1,0)$};
	\node[right= 0.2 of 12110] (21101) {$(2,1,1,0,1)$};
	\node[right= 0.2 of 21101] (12101) {$(1,2,1,0,1)$};
	\node[right= 0.2 of 12101] (21011) {$(2,1,0,1,1)$};
	\node[right= 0.2 of 21011] (12011) {$(1,2,0,1,1)$};
	\node[right= 0.2 of 12011] (02111) {$(2,0,1,1,1)$};
	\node[above = 1.9 of 20111] (21111) {$(2,1,1,1,1)$}; 
	\node[above right = 1.9 and -0.9 of 21110] (22110) {$(2,2,1,1,0)$}; 
	\node[above right = 1.9 and -0.9 of 21101] (22101) {$(2,2,1,0,1)$}; 
	\node[above right = 1.9 and -0.9 of 21011] (22011) {$(2,2,0,1,1)$}; 
	\node[above = 1.9 of 02111] (12111) {$(1,2,1,1,1)$}; 
	\node[above = 1.9 of 22101] (22111) {$(2,2,1,1,1)$}; 
	\node[below = 3.9  of 22101] (0hat) {$\widehat{0}$}; 
\draw (20111) to (21111);
\draw (11111.north) to (21111);
\draw (21110) to (21111);
\draw (21101.north) to (21111);
\draw (21011.north) to (21111);
\draw (11111.north) to (12111);
\draw (12110.north) to (12111);
\draw (12101.north) to (12111);
\draw (12011) to (12111);
\draw (02111) to (12111);
\draw (21110.north) to (22110);
\draw (12110.north) to (22110);
\draw (21101.north) to (22101);
\draw (12101.north) to (22101);
\draw (21011.north) to (22011);
\draw (12011.north) to (22011);
\draw (21111) to (22111);
\draw (22110) to (22111);
\draw (22101) to (22111);
\draw (22011) to (22111);
\draw (12111) to (22111);
\draw (0hat) to (20111.south);
\draw (0hat) to (11111.south);
\draw (0hat) to (21110.south);
\draw (0hat) to (12110);
\draw (0hat) to (21101);
\draw (0hat) to (12101);
\draw (0hat) to (21011);
\draw (0hat) to (12011.south);
\draw (0hat) to (02111.south);
\end{tikzpicture}}.
 \]
 
 By Theorem~\ref{thm:main}, we may compute the Grothendieck polynomial $\mathfrak{G}_{143562}$ by 
 \begin{itemize}
 	\item  removing the minimum element $\hat{0}$ of $P_{143562}$,
 	\item  labelling each remaining minimal element by $1$, 
 	\item  and labelling all other element in such a way that, for each $p \in P_{143562}$, the sum of the labels of all $q \leq p$ equals $1$. 
 \end{itemize}
 In this case, our labelled truncated poset diagram becomes
  \[\scalebox{1}
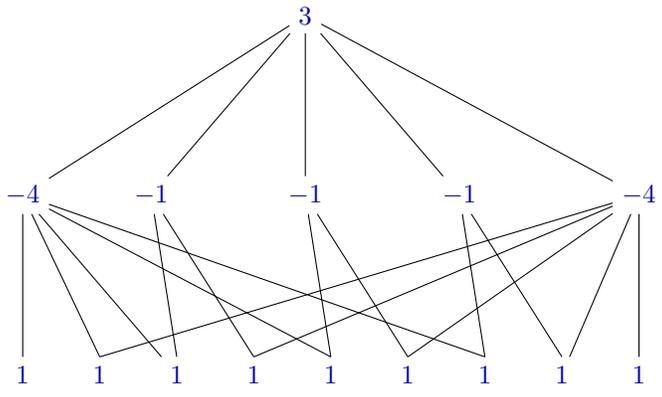
{
 \begin{tikzpicture}
	\node (20111) at (0,0) {$\textcolor{darkblue}{1}$};
	\node[right= 0.6 of 20111] (11111) {$\textcolor{darkblue}{1}$};
	\node[right= 0.6 of 11111] (21110) {$\textcolor{darkblue}{1}$};
	\node[right= 0.6 of 21110] (12110) {$\textcolor{darkblue}{1}$};
	\node[right= 0.6 of 12110] (21101) {$\textcolor{darkblue}{1}$};
	\node[right= 0.6 of 21101] (12101) {$\textcolor{darkblue}{1}$};
	\node[right= 0.6 of 12101] (21011) {$\textcolor{darkblue}{1}$};
	\node[right= 0.6 of 21011] (12011) {$\textcolor{darkblue}{1}$};
	\node[right= 0.6 of 12011] (02111) {$\textcolor{darkblue}{1}$};
	\node[above = 1.9 of 20111] (21111) {$\textcolor{darkblue}{-4}$}; 
	\node[above right = 1.9 and -0.9 of 21110] (22110) {$\textcolor{darkblue}{-1}$}; 
	\node[above right = 1.9 and -0.9 of 21101] (22101) {$\textcolor{darkblue}{-1}$}; 
	\node[above right = 1.9 and -0.9 of 21011] (22011) {$\textcolor{darkblue}{-1}$}; 
	\node[above = 1.9 of 02111] (12111) {$\textcolor{darkblue}{-4}$}; 
	\node[above = 1.9 of 22101] (22111) {$\textcolor{darkblue}{3}$}; 
\draw (20111) to (21111);
\draw (11111.north) to (21111);
\draw (21110) to (21111);
\draw (21101.north) to (21111);
\draw (21011.north) to (21111);
\draw (11111.north) to (12111);
\draw (12110.north) to (12111);
\draw (12101.north) to (12111);
\draw (12011) to (12111);
\draw (02111) to (12111);
\draw (21110.north) to (22110);
\draw (12110.north) to (22110);
\draw (21101.north) to (22101);
\draw (12101.north) to (22101);
\draw (21011.north) to (22011);
\draw (12011.north) to (22011);
\draw (21111) to (22111);
\draw (22110) to (22111);
\draw (22101) to (22111);
\draw (22011) to (22111);
\draw (12111) to (22111);
\end{tikzpicture}},
 \]
 so we obtain that 
 \begin{align*}
 	\mathfrak{G}_{143562} &= \mathfrak{S}_{143562} - 4\xx^{(2,1,1,1,1)} - \xx^{(2,2,1,1,0)} - \xx^{(2,2,1,0,1)} - \xx^{(2,2,0,1,1)} - 4\xx^{(1,2,1,1,1)} + 3 \xx^{(2,2,1,1,1)}.
 \end{align*}
This calculation may be verified by comparison with other combinatorial formulas for the Grothendieck polynomial $\mathfrak{G}_{143562}$. For another example (with $w=351624$) of the application of Theorem~\ref{thm:main}, see \cite[Example~4.2 and Figure~7]{MSSD}.

\section*{Acknowledgements}
We thank Jenna Rajchgot for pointing us to helpful references. OP is grateful for conversations about \cite{MSSD21} with Zach Hamaker and Anna Weigandt.

\providecommand{\bysame}{\leavevmode\hbox to3em{\hrulefill}\thinspace}
\providecommand{\MR}{\relax\ifhmode\unskip\space\fi MR }
\providecommand{\MRhref}[2]{%
  \href{http://www.ams.org/mathscinet-getitem?mr=#1}{#2}
}
\providecommand{\href}[2]{#2}

\end{document}